\documentclass[fleqn,10pt]{wlscirep}
\usepackage[utf8]{inputenc}
\usepackage[T1]{fontenc}
\title{Unbalanced Regularized Optimal Mass Transport with Applications to Fluid Flows in the Brain}

\author[1,*]{Xinan Chen}
\author[2]{Helene Benveniste}
\author[3]{Allen R. Tannenbaum}
\affil[1]{Memorial Sloan Kettering Cancer Center, Department of Medical Physics, New York, 10065, USA}
\affil[2]{Yale School of Medicine, Department of Anesthesiology, New Haven, 06510, USA}
\affil[3]{Stony Brook University, Departments of Computer Science and Applied Mathematics \& Statistics, Stony Brook, 11794, USA}

\affil[*]{corresponding author; chenx7@mskcc.org}

\begin{abstract}
As a generalization of the optimal mass transport (OMT) approach of Benamou and Brenier's, the regularized optimal mass transport (rOMT) formulates a transport problem from an initial mass configuration to another with the optimality defined by the total kinetic energy, but subject to an advection-diffusion constraint equation. Both rOMT and the Benamou and Brenier's formulation require the total initial and final masses to be equal; mass is preserved during the entire transport process. However, for many applications, e.g., in dynamic image tracking, this constraint is rarely if ever satisfied. Therefore, we propose to employ an unbalanced version of rOMT to remove this constraint together with a detailed numerical solution procedure and applications to analyzing fluid flows in the brain.
\end{abstract}
\begin{document}

\flushbottom
\maketitle
%
%
\thispagestyle{empty}


\section*{Introduction}

The optimal mass transport (OMT) problem is concerned with finding a transport mapping from an initial mass distribution to a final one, with optimality defined relative to a given cost function \cite{GM,LK,villani1,villani2}. A reformulation of the OMT problem in a fluid dynamical framework was proposed by Benamou and Brenier \cite{benamou2000computational} using the $L^2$ distance as the basis for the cost function, which may be written as follows in the special case of interest to us in the present work. 

Given an initial mass distribution function $\rho_0(x)\geqslant0$ and a final one $\rho_1(x)\geqslant0$ both defined on a bounded region $\Omega\subseteq\mathbb{R}^3$ and with the same total mass (i.e., $\int_{\Omega}\rho_0(x)\,\mathrm{d}x = \int_{\Omega}\rho_1(x)\,\mathrm{d}x$), the dynamic OMT problem by Benamou and Brenier \cite{benamou2000computational} aims to solve
\begin{subequations}\label{eq:BB}
\begin{align}
\underset{\rho,v}{\text{min}}\quad & \int_0^T\int_{\Omega}\left\lVert v(t,x)\right\rVert^2\rho(t,x)\,\mathrm{d}x\,\mathrm{d}t\label{eq:BB_energy}\\
\text{subject to}\quad
& \frac{\partial\rho}{\partial t} + \nabla\cdot(\rho v) = 0,\label{eq:BB_conti}\\
& \rho(0,x) = \rho_0(x), \quad\rho(T,x) = \rho_1(x) \label{eq:BB_bound}
\end{align}
\end{subequations}
where a temporal dimension $t\in[0,T]$ is added to the transport process. In the above formulation, $\rho(t,x)$ is the dynamic density function, and  $v(t,x)$ is the dynamic velocity field defining the fluid flows from $\rho_0$ to $\rho_1$. Equation \eqref{eq:BB_conti} is called the \textit{continuity equation} in fluid dynamics, and characterizes the advective transport of a conserved quantity in bulk flows. The cost function \eqref{eq:BB_energy} is the total kinetic energy of the transport process. The square root of the achieved minimum of \eqref{eq:BB_energy}, if it exists, is called the $L^2$-Wasserstein distance between $\rho_0$ and $\rho_1$.

As an extension of model \eqref{eq:BB_energy}-\eqref{eq:BB_bound}, a regularized version of OMT (rOMT) has been developed and applied in various places in which one includes a diffusion motion into the transport; see relavent work\cite{cuturi2013sinkhorn,pavon,koundal2020} and the many references therein. More precisely, a diffusion term is added into the continuity equation \eqref{eq:BB_conti} to make it
\begin{equation} \label{eq:Diff}
\frac{\partial\rho}{\partial t} + \nabla\cdot(\rho v) = \sigma\Delta\rho
\end{equation}
where the constant $\sigma>0$ is the diffusion coefficient. Equation~\eqref{eq:Diff} is thus an \textit{advection-diffusion equation} in fluid dynamics. The rOMT model has been proven useful for a number of important tracking problems in computational fluid dynamics, such as in quantifying and visualizing the movement of solutes in the brain on dynamic contrast enhanced magnetic resonance imaging (DCE-MRI)\cite{koundal2020,chen2022,chen2022visualizing1,ROBERT2023764,ozturk2023continuous}.

As is well-known, both the OMT and rOMT models must satisfy the total mass conservation constraint, namely $\int_{\Omega}\rho_0(x)\,\mathrm{d}x = \int_{\Omega}\rho_1(x)\,\mathrm{d}x$. In this case, we usually call the problem as \textit{balanced} to refer to the conservation of total mass. Indeed as formulated in \eqref{eq:Diff}, neither advection nor diffusion will change the total mass locally or globally in a given region $\Omega$. However, for applications in dynamic imagery in which either OMT or rOMT is employed as an optical flow tracking method, this is almost never the case. For example, in DCE-MRI data where a gadolinium-based tracer is injected and delivered into the body, an early climbing period of the total image intensity is usually observed, since it takes time for the tracer to reach and fill the region of interest. Under these circumstances, if we assume that the intensity of image signal is proportional to the density and mass in the aforementioned two models, the total mass conservation law is no longer satisfied, and thus the OMT and rOMT models cannot be directly applied. Notably, analyzing the initial accumulating stage of tracers may help uncover interesting and physiologically relevant transport patterns in the brain, and therefore a new model is necessary. 

In the present work,  we propose an \textit{unbalanced regularized OMT} (urOMT) model for applications to the DCE-MRI fluid flow data, where an independent variable and its indicator function are added as an ``invisible" sink or source of mass. The urOMT problem is formulated as follows:
\begin{subequations}\label{eq:uromt_dp}
\begin{align}
\underset{\rho,v,r}{\text{min}}\quad & \int_0^T\int_{\Omega}
\begin{aligned}[t]\label{eq:uromt_dp1}
& \big(\left\lVert v(t,x)\right\rVert^2\rho(t,x) \\
& + \alpha \chi(t,x)r(t,x)^2\rho(t,x)\big)\,\mathrm{d}x\,\mathrm{d}t\\
\end{aligned}\\
\text{subject to}\quad
& \frac{\partial\rho}{\partial t} + \nabla\cdot(\rho v) = \sigma\Delta\rho + \chi\rho r, \label{eq:uromt_conti}\\
& \rho(0,x) = \rho_0(x), \quad\rho(T,x) = \rho_1(x) \label{eq:uromt_dp3}
\end{align}
\end{subequations}
where $r(t,x)$ is the relative source variable, $\chi(t,x)$ is the given indicator function of $r(t,x)$ which takes values either 0 or 1 to constrain $r$ to a certain spatial and temporal location, and $\alpha>0$ is the weighting parameter of the source term in the cost function. This model takes inputs $\rho_0(x)$, $\rho_1(x)$ and $\chi(t,x)$, and solves for the optimal $\rho(t,x)$, $v(t,x)$ and $r(t,x)$. The added second term in the cost function \eqref{eq:uromt_dp1} is called the \textit{Fisher-Rao} term which arises from the Fisher-Rao metric in information geometry. Our model therefore can be viewed as the interpolation between the $L^2$-Wasserstein and the Fisher-Rao metrics\cite{chizat2018interpolating,CheGeoTan16b}. 

The partial differential equation \eqref{eq:uromt_conti} indicates that there are three types of physical phenomena taking place in the dynamic system, advection ($\nabla\cdot(\rho v)$), diffusion ($\sigma\Delta\rho$) and mass creation/destruction ($\chi\rho r$). We should note that we prefer to use the relative source $r$ which controls the rate of mass gain ($r>0$) and loss ($r<0$), rather than the source $s=\rho r$, as the unbalanced variable in our urOMT formulation, since the resulting expressions are cleaner and easier to interpret. Even though $r$ plays a role as a sink of mass when $r<0$ and a role as a source of mass when $r>0$, in our work we refer to $r$ as the relative ``source'' to broadly refer to its ability to generate unbalanced mass instantaneously. One can imagine that $r$ is the source of both positive and negative mass. The main point however, is that we no longer require the total mass conservation condition for the input images $\rho_0$ and $\rho_1$.

The unbalanced OMT problem has been studied both theoretically \cite{chizat2018unbalanced,CheGeoTan16b} and with various applications in meteorology \cite{benamou}, shape modification \cite{chizat2018scaling}, image registration \cite{feydy2017optimal}, image deformation \cite{maas2015generalized,sejourne2022unbalanced}, tumor growth modeling \cite{lombardi2015eulerian,galloute2017unbalanced}, population modeling \cite{yang2018scalable} and dynamical tracking \cite{lee2019parallel}, etc.

Our work presented here should be considered as an extension of the rOMT algorithm \cite{chen2022visualizing1} since both share the similar numerical structure and setup. We are specifically interested in applying the urOMT method into DCE-MRI studies to quantify the fluid flows in the rat brain which is our motivation of introducing an unbalanced term to analyze unbalanced data. Although there has been a good amount of work in unbalanced OMT as mentioned above, to the best of our knowledge, this is the first work to incorporate the unbalanced regularized OMT for the quantification of dynamic fluid flows using DCE-MRI imaging.

Briefly summarizing the present paper, in the Section "Numerical Method", we give the detailed numerical method for solving the urOMT problem \eqref{eq:uromt_dp1}-\eqref{eq:uromt_dp3}, and in the Section "Results", we test our urOMT model and show results of applications to studying dynamic fluid flows in both synthetic data and DCE-MRI rat brain data. In the Section "Discussion", we further discuss and summarize the urOMT method. Some relevant efforts and potential future work are also provided. Lastly in the Section "Conclusion", we conclude our present work.



\section*{Numerical Method}\label{sec:num}
In this section, we elaborate on the numerical method developed for the urOMT problem \eqref{eq:uromt_dp1}-\eqref{eq:uromt_dp3} on 3D images, especially on DCE-MRI-based images in which the signal intensity levels are reflecting the concentration of the Gadolinium based tracer, and are therefore proportional to the mass/density in our model. Note that this method can be easily adapted to images in any dimension with minor changes. The numerical method used in this work is largely inspired by and based on the previous work\cite{chen2022visualizing1}.

\subsection*{Model}
Given a pair of 3D images, $\rho_0^{img}(x)$ and $\rho_1^{img}(x)$, and an indicator $\chi(x,t)$, in avoidance of the over-matching of the image noise, we consider posing a free end-point condition, so that we can remove the end-point constraint $\rho(T,x) = \rho_1^{img}(x)$. Another fitting term is consequently added into the cost function. The numerical model we solve is therefore written as:
\begin{subequations}\label{eq:uromt_dp2}
\begin{align}
\underset{v,r}{\text{min}}\quad &
\begin{aligned}[t]\label{eq:uromt2_func}
&\int_0^T\int_{\Omega}\rho(t,x)\big(\left\lVert v(t,x)\right\rVert^2 \\
&+ \alpha \chi(t,x)r(t,x)^2\big)\,\mathrm{d}x\,\mathrm{d}t \\
& + \beta\int_{\Omega}(\rho(T,x)-\rho_1^{img}(x))^2\,\mathrm{d}x\\
\end{aligned}
\\
\text{subject to}\quad
& \frac{\partial\rho}{\partial t} + \nabla\cdot(\rho v) = \sigma\Delta\rho + \chi\rho r, \label{eq:uromt2_conti}\\
& \rho(0,x) = \rho_0^{img}(x)\label{eq:uromt2_bound}
\end{align}
\end{subequations}
where $\beta>0$ is the weighting parameter for the fitting term in the cost function. The dynamic density function $\rho(t,x)$ can be explicitly derived starting at $\rho_0^{img}$ and following equation \eqref{eq:uromt2_conti} with a velocity field $v$ and a relative source $r$ with its indicator $\chi$, so $\rho$ is removed from the optimized variables. Basically, with this setup $\rho$ becomes a {\em state variable}. We call this model \eqref{eq:uromt2_func}-\eqref{eq:uromt2_bound} the \textit{unbalanced regularized OMT with free end-point} (free-urOMT).

\subsection*{Discretization}
Since 3D images are typically defined on cubical domains, we divide the cubical space $\Omega$ into a cell-centered grid size $n_1\times n_2\times n_3$ with uniform spacing $\Delta x$, $\Delta y$ and $\Delta z$ in $x$, $y$ and $z$-direction, respectively. Let $n=n_1n_2n_3$ be the total number of voxels. The time interval $[0,T]$ is partitioned into $m$ equal sub-intervals with length $\Delta t=\frac{T}{m}$. Then we have $m+1$ discrete time steps $t_i = i\Delta t$ for $i = 0,\cdots,m$. 

We use a bold font to denote a flattened vector discretized from its corresponding continuous function onto the cell-centered grid defined above. Therefore, the given initial and final images $\rho_0^{img}(x)$ and $\rho_1^{img}(x)$ is discretized into vectors $\pmb{\rho_0^{img}}$ and $\pmb{\rho_1^{img}}$. The density function $\rho(t,x)$ is discretized into $\pmb{\rho_i}$ for $i=0,\cdots,m$, each denoting the mass distribution at $t_i$, and where $\pmb{\rho_0} = \pmb{\rho_0^{img}}$ denotes the given initial image. The velocity field $v(t,x)$, relative source $r(t,x)$ and its indicator $\chi(t,x)$ may also be discretized into $\pmb{v_i}$, $\pmb{r_i}$ and $\pmb{\chi_i}$ for $i=0,\cdots,m-1$, each denoting the velocity field, relative source, and the indicator transforming $\pmb{\rho_{i}}$ to $\pmb{\rho_{i+1}}$, respectively. We further denote $\pmb{\rho} = [\pmb{\rho_1};\cdots;\pmb{\rho_m}]$, $\pmb{v} = [\pmb{v_0};\cdots;\pmb{v_{m-1}}]$, $\pmb{r} = [\pmb{r_0};\cdots;\pmb{r_{m-1}}]$ and $\pmb{\chi} = [\pmb{\chi_0};\cdots;\pmb{\chi_{m-1}}]$.

So far, we have defined the discretized variables on the space and time grids. The cost function equation~\eqref{eq:uromt2_func} may therefore be approximated by
\begin{equation}
\Gamma(\pmb{v},\pmb{r}) = \Gamma_1(\pmb{v},\pmb{r}) + \alpha\Gamma_2(\pmb{v},\pmb{r}) + \beta\Gamma_3(\pmb{v},\pmb{r})
\end{equation}
where
\begin{subequations}
\begin{align}
&\Gamma_1(\pmb{v},\pmb{r}) = (\Delta t\Delta x\Delta y\Delta z)\pmb{\rho}^{T}(I_m\otimes[I_n|I_n|I_n])(\pmb{v}\odot \pmb{v}), \label{eq:gamma1}\\
&\Gamma_2(\pmb{v},\pmb{r}) = (\Delta t\Delta x\Delta y\Delta z)\pmb{\rho}^{T}(\pmb{r}\odot \pmb{r}\odot\pmb{\chi}),\\
&\Gamma_3(\pmb{v},\pmb{r}) = (\Delta x\Delta y\Delta z)\|\pmb{\rho_m} - \pmb{\rho_1^{img}}\|^2. \label{eq:gamma3}
\end{align}
\end{subequations}
Here $\otimes$ denotes the Kronecker tensor product, and $\odot$ denotes the Hadamard product. Further, $[\cdot|\cdot]$ represents the block matrix, and $\|\cdot\|$ means taking the $L^2$ norm of a vector. $I_k$ is the $k$-dimensional identity matrix for $k\in\mathbb{N}^+$.

\subsection*{Solving the Partial Differential Equation}

Next, we deal with the partial differential equation \eqref{eq:uromt2_conti}. We place ghost points outside of the boundary and employ the Neumann boundary condition such that the derivative across the boundary is always 0. The Laplacian operator $\Delta$ may be approximated with a matrix $Q$ under the aforementioned numerical grid and boundary condition. Similar to the technique employed in the previous work\cite{chen2022visualizing1}, we use the operator-splitting method to numerically solve the equation but in this work we divide the equation into three steps. To be precise, at each time step from $t_i$ to $t_{i+1}$ for $i = 0,\cdots,m-1$, we divide the whole process \eqref{eq:uromt2_conti} into first, mass gain/loss: $\frac{\partial\rho}{\partial t} = \chi\rho r$; second, advection: $\frac{\partial\rho}{\partial t} + \nabla\cdot(\rho v) = 0$; and third, diffusion: $\frac{\partial\rho}{\partial t} = \sigma\Delta\rho$, in total three steps, and then integrate them together.

For the first mass gain/loss step, given an initial condition $\rho(t_i, x) = \rho_i(x)$, the equation may be discretized into
\begin{subequations}
\begin{align}
& \frac{1}{\Delta t}(\pmb{\rho_i^{src}}-\pmb{\rho_i}) = \pmb{\rho_i} \odot \pmb{r_i}\odot\pmb{\chi_i}\\
\Rightarrow\quad & \pmb{\rho_i^{src}} = (\pmb{1_n} + \Delta t\cdot \pmb{r_i}\odot\pmb{\chi_i})\odot \pmb{\rho_i} \label{eq:src}
\end{align}
\end{subequations}
where $\pmb{1_n}$ is a vector of length $n$ consisting of 1's. The second advection step with an initial condition $\rho(t_i, x) = \rho_i^{src}(x)$ can be discretized and solved with
\begin{equation}\label{eq:adv}
\pmb{\rho_i^{adv}} = S(\pmb{v_i})\pmb{\rho_i^{src}},
\end{equation}
where $S(\pmb{v_i})$ is the averaging matrix with respect to $\pmb{v_i}$ using the particle-in-cell method which redistributes the transported mass to its nearest neighbors by a certain ratio. The $(j,k)$ entry of $S(\pmb{v_i})$ is the ratio of mass allocated from the old location $k$ to the new location $j$. Multiplying $S(\pmb{v_i})$ by a vector which represents a mass distribution, we can derive a new mass distribution transported by the velocity $\pmb{v_i}$ under the pre-defined numerical grid. The third diffusion step with an initial condition $\rho(t_i, x) = \rho_i^{adv}(x)$ employs the Euler Backwards scheme in the following manner:
\begin{subequations}
\begin{align}
& \frac{1}{\Delta t}(\pmb{\rho_{i+1}}-\pmb{\rho_i^{adv}}) = \sigma Q\pmb{\rho_{i+1}}\\
\Rightarrow\quad & (I_n-\sigma \Delta t\cdot Q)\pmb{\rho_{i+1}} = \pmb{\rho_i^{adv}}. \label{eq:diff}
\end{align}
\end{subequations}
Combining all three steps \eqref{eq:src}, \eqref{eq:adv} and \eqref{eq:diff}, we have
\begin{equation}\label{eq:eq1}
(I_n-\sigma \Delta t\cdot Q)\pmb{\rho_{i+1}} = S(\pmb{v_i})(\pmb{1_n} + \Delta t\cdot \pmb{r_i}\odot\pmb{\chi_i})\odot \pmb{\rho_i}.
\end{equation}

\begin{figure*}
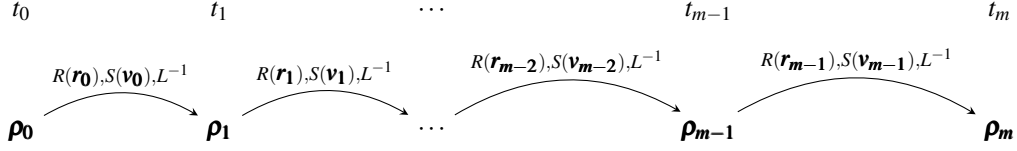

\begin{center}
\begin{codi}[l/.style={bend left}, r/.style={bend right}]
  \obj {t_0 & [3 em] t_1 & [3.5 em] |(pb)| \cdots & [6 em] t_{m-1} & [6.5 em] t_{m} \\
           |(mu0)| \pmb{\rho_0} & |(mu1)| \pmb{\rho_1} & |(pbs)| \cdots & |(mum-1)| \pmb{\rho_{m-1}} & |(mum)| \pmb{\rho_{m}} \\ };

  \mor  mu0 "\large{R(\pmb{r_0}),S(\pmb{v_0}),L^{-1}}":l,->   mu1 "\large{R(\pmb{r_1}),S(\pmb{v_1}),L^{-1}}":l,-> pbs "\large{R(\pmb{r_{m-2}}),S(\pmb{v_{m-2}}),L^{-1}}":l,->  mum-1 "\large{R(\pmb{r_{m-1}}),S(\pmb{v_{m-1}}),L^{-1}}":l,->  mum;  
\end{codi}
\end{center}
\caption{Numerical Pipeline of urOMT. From $t_i$ to $t_{i+1}$ for $i = 0,\cdots,m-1$, the interpolated image $\pmb{\rho_i}$ is firstly added with mass by applying matrix $R(\pmb{r_i})$, and is secondly advected via the velocity field $\pmb{v_i}$ by applying the averaging matrix $S(\pmb{v_i})$ and is lastly diffused by applying matrix $L^{-1}$.} \label{fig:uromt.pipeline}
\end{figure*}

If we denote
\begin{subequations}
\begin{align}
& L \triangleq I_n-\sigma \Delta t\cdot Q,\\
& R(\pmb{r_i}) \triangleq I_n+\Delta t\cdot \text{diag}(\pmb{r_i}\odot\pmb{\chi_i}),
\end{align}
\end{subequations}
where $\text{diag}(\cdot)$ is the operator turning a vector into a diagonal matrix. Equation~\eqref{eq:eq1} may be re-written as 
\begin{equation}
\pmb{\rho_{i+1}} = L^{-1}S(\pmb{v_i})R(\pmb{r_i})\pmb{\rho_i}
\end{equation}
for $i=0,\cdots,m-1$ where $S(\pmb{v_i})$ and $R(\pmb{r_i})$ are $n$-dimensional matrix with respect to $\pmb{v_i}$ and $\pmb{r_i}$, respectively. See Fig \ref{fig:uromt.pipeline} for the pipeline of the numerical dynamics.

In conclusion, we have discretized the free-urOMT problem \eqref{eq:uromt2_func}-\eqref{eq:uromt2_bound} as follows:
\begin{subequations}\label{eq:uromt_dp2_disc}
\begin{align}
\underset{\pmb{v},\pmb{r}}{\text{min}}\quad & \Gamma(\pmb{v},\pmb{r})=\Gamma_1(\pmb{v},\pmb{r})+\alpha\Gamma_2(\pmb{v},\pmb{r})+\beta\Gamma_3(\pmb{v},\pmb{r}) \label{eq:uromt_dp2_disc_func}\\
\text{subject to}\quad
& \pmb{\rho_{i+1}} = L^{-1}S(\pmb{v_i})R(\pmb{r_i})\pmb{\rho_i} \text{, for } i=0,\cdots,m-1,\label{eq:chain}\\
& \pmb{\rho_0} = \pmb{\rho_0^{img}}, \label{eq:uromt_dp2_disc_bound}
\end{align}
\end{subequations}
where $\Gamma_1, \Gamma_2$ and $\Gamma_3$ are explicitly given in \eqref{eq:gamma1}-\eqref{eq:gamma3}.

\subsection*{Computing the Gradient and the Hessian}\label{sec:gradient}
Notice that in the discrete problem \eqref{eq:uromt_dp2_disc_func}-\eqref{eq:uromt_dp2_disc_bound}, we can equivalently re-write the constraint \eqref{eq:chain} as
\begin{equation}\label{eq:chainsss}
\pmb{\rho_{i}} = L^{-1}S(\pmb{v_{i-1}})R(\pmb{r_{i-1}})L^{-1}S(\pmb{v_{i-2}})R(\pmb{r_{i-2}})\cdots L^{-1}S(\pmb{v_{0}})R(\pmb{r_{0}})\pmb{\rho_0}
\end{equation}
for $i=1,\cdots,m$, which explicitly gives the expression of mass distributions at all time steps from a given initial mass distribution $\pmb{\rho_0}$, a velocity field $\pmb{v}$ and a relative source $\pmb{r}$. One can prove that $S(\pmb{v_{i}})$ is linear to $\pmb{v_i}$ and $R(\pmb{v_{i}})$ is linear to $\pmb{r_i}$ for $i=0,\cdots,m-1$. Then according to equations \eqref{eq:chainsss} and \eqref{eq:uromt_dp2_disc_bound}, $\pmb{\rho_m}$ and $\pmb{\rho}=[\pmb{\rho_1};\cdots;\pmb{\rho_m}]$ in $\Gamma(\pmb{v},\pmb{r})$ can both be explicitly written in linear to $\pmb{v}$ and $\pmb{r}$. Therefore, the optimization problem \eqref{eq:uromt_dp2_disc_func}-\eqref{eq:uromt_dp2_disc_bound} can indeed be viewed as an unconstrained minimization problem. Observe that in the cost function, $\Gamma_1$ is linear to $\pmb{r}$ and its component $\pmb{v}\odot\pmb{v}$ is quadratic to $\pmb{v}$; $\Gamma_2$ is linear to $\pmb{v}$ and its component $\pmb{r}\odot\pmb{r}$ is quadratic to $\pmb{r}$; $\Gamma_3$ is quadratic to both $\pmb{v}$ and $\pmb{r}$. It is then natural to employ the Gauss-Newton method to optimize on the problem, which involves computing the gradient and the Hessian matrix of $\Gamma$ with respect to variables $\pmb{v}$ and $\pmb{r}$. 

Next, we focus on calculating the gradient of $\Gamma$: 
\begin{equation}
g\triangleq[g_{\pmb{v}};g_{\pmb{r}}] \text{, where } g_{\pmb{v}}=\frac{\partial\Gamma}{\partial \pmb{v}}, g_{\pmb{r}}=\frac{\partial\Gamma}{\partial \pmb{r}}
\end{equation}
and the Hessian matrix of $\Gamma$: 
\begin{equation}
H\triangleq
\begin{pmatrix}
    H_{11} & H_{12} \\
    H_{21} & H_{22} 
\end{pmatrix}
\text{, where }
H_{11} = \frac{\partial^2\Gamma}{\partial \pmb{v}^2},\quad
H_{12} = \frac{\partial^2\Gamma}{\partial \pmb{v}\partial \pmb{r}},\quad
H_{21} = \frac{\partial^2\Gamma}{\partial \pmb{r}\partial \pmb{v}},\quad
H_{22} = \frac{\partial^2\Gamma}{\partial \pmb{r}^2}.
\end{equation}

Equation~\eqref{eq:chainsss} shows that $\pmb{\rho_{k}}$ is determined by $\pmb{v_0},\cdots,\pmb{v_{k-1}}, \pmb{r_0},\cdots,\pmb{r_{k-1}}$ and is thus independent of $\pmb{v_j}$ and $\pmb{r_j}$ for $j\geqslant k$. Defining 
\begin{equation}
J_{\pmb{v_j}}^k\triangleq\frac{\partial\pmb{\rho_k}}{\partial \pmb{v_j}},\quad 
J_{\pmb{r_j}}^k\triangleq\frac{\partial\pmb{\rho_k}}{\partial \pmb{r_j}}
\end{equation}
for $k = 1,\cdots,m, j = 0,\cdots,m-1$, then $J_{\pmb{v_j}}^k = 0$ and $J_{\pmb{r_j}}^k = 0$ always hold for $j\geqslant k$. If we further denote
\begin{equation}
J_{\pmb{v}} \triangleq \frac{\partial\pmb{\rho}}{\partial \pmb{v}} = (J_{\pmb{v_j}}^k)_{k,j},\quad J_{\pmb{r}} \triangleq \frac{\partial\pmb{\rho}}{\partial \pmb{r}} = (J_{\pmb{r_j}}^k)_{k,j},
\end{equation}
then $J_{\pmb{v}}$ and $J_{\pmb{r}}$ are lower-triangular block matrices of the form

\begin{equation}
    J_{\pmb{v}} =\begin{pmatrix}
        J_{\pmb{v_0}}^1 &  & & \\
        J_{\pmb{v_0}}^2 & J_{\pmb{v_1}}^2 & & \\
        \vdots & \vdots & \ddots &\\
        J_{\pmb{v_0}}^m & J_{\pmb{v_1}}^m & \cdots & J_{\pmb{v_{m-1}}}^m\\
        \end{pmatrix}
    \triangleq
    \begin{pmatrix}
        J^1_{\pmb{v}} \\
        J^2_{\pmb{v}} \\
        \vdots\\
        J^m_{\pmb{v}} \\
        \end{pmatrix}
\text{\quad and \quad}
    J_{\pmb{r}} =\begin{pmatrix}
        J_{\pmb{r_0}}^1 &  & & \\
        J_{\pmb{r_0}}^2 & J_{\pmb{r_1}}^2 & & \\
        \vdots & \vdots & \ddots &\\
        J_{\pmb{r_0}}^m & J_{\pmb{r_1}}^m & \cdots & J_{\pmb{r_{m-1}}}^m\\
        \end{pmatrix}
    \triangleq
    \begin{pmatrix}
        J^1_{\pmb{r}} \\
        J^2_{\pmb{r}} \\
        \vdots\\
        J^m_{\pmb{r}} \\
        \end{pmatrix}
\end{equation}
where $J_{\pmb{v}}^k=\big[J_{\pmb{v_0}}^k|J_{\pmb{v_1}}^k|\cdots |J_{\pmb{v_{k-1}}}^k\big]$ denotes the row block of $J_{\pmb{v}}$, and $J_{\pmb{r}}^k$ for that of $J_{\pmb{r}}$ for $k=1,\cdots,m$. With the notations defined above, then for the gradients we have
\begin{subequations}\label{eq:gradients1}
\begin{align}
g_{\pmb{v}} & = \frac{\partial\Gamma_1}{\partial \pmb{v}} + \alpha\frac{\partial\Gamma_2}{\partial \pmb{v}} + \beta\frac{\partial\Gamma_3}{\partial \pmb{v}}\\
& = 
\begin{aligned}[t]\label{eq:eq3}
& (\Delta t \Delta x \Delta y \Delta z)\left(2(M\text{diag}(\pmb{v}))^T\pmb{\rho} + J_{\pmb{v}}^TM(\pmb{v}\odot \pmb{v}) \right)\\
& + \alpha (\Delta t \Delta x \Delta y \Delta z)J_{\pmb{v}}^T(\pmb{r}\odot \pmb{r}\odot\pmb{\chi})\\
& + 2\beta(\Delta x \Delta y \Delta z)(J_{\pmb{v}}^m)^T\big(\pmb{\rho_m} - \pmb{\rho_1^{img}}\big)
\end{aligned}
\end{align}
\end{subequations}
and 
\begin{subequations}\label{eq:gradients2}
\begin{align}
g_{\pmb{r}} & = \frac{\partial\Gamma_1}{\partial \pmb{r}} + \alpha\frac{\partial\Gamma_2}{\partial \pmb{r}} + \beta\frac{\partial\Gamma_3}{\partial \pmb{r}}\\
& = 
\begin{aligned}[t]
& (\Delta t \Delta x \Delta y \Delta z)J_{\pmb{r}}^TM(\pmb{v}\odot \pmb{v})\\
& + \alpha (\Delta t \Delta x \Delta y \Delta z)\left(2\text{diag}(\pmb{r}\odot\pmb{\chi})\pmb{\rho}+J_{\pmb{r}}^T(\pmb{r}\odot \pmb{r}\odot\pmb{\chi})\right)\\
& + 2\beta(\Delta x \Delta y \Delta z)(J_{\pmb{r}}^m)^T\big(\pmb{\rho_m} - \pmb{\rho_1^{img}}\big)
\end{aligned}
\end{align}
\end{subequations}
where we use $M \triangleq I_m\otimes[I_n|I_n|I_n]$ for simplicity of notation. For the Hessian matrix $H$, we use a function handle which is a MATLAB data type in anticipation of solving the linear system $H\pmb{x}=-g$ in the next step. To be specific, we compute vector $H\pmb{x}$ where $\pmb{x}=[\pmb{x_v};\pmb{x_r}]$ and $\pmb{x_v}\in\mathbb{R}^{3mn}, \pmb{x_r}\in\mathbb{R}^{mn}$ rather than matrix $H$. To satisfy the symmetry of the Hessian matrix and to omit some complex second-order terms, we have the approximations
\begin{subequations}\label{eq:hessians}
\begin{align}
H_{11}\pmb{x_v} & \approx 2(\Delta t \Delta x \Delta y \Delta z)\text{diag}(\pmb{\rho}^{T}M) + 2\beta(\Delta x \Delta y \Delta z)(J_{\pmb{v}}^m)^TJ_{\pmb{v}}^m\pmb{x_v},\\
H_{22}\pmb{x_r} & \approx 2\alpha(\Delta t \Delta x \Delta y \Delta z)\text{diag}(\pmb{\rho}\odot\pmb{\chi}) + 2\beta(\Delta x \Delta y \Delta z)(J_{\pmb{r}}^m)^TJ_{\pmb{r}}^m\pmb{x_r},\\
H_{12}\pmb{x_r} & \approx 2\beta(\Delta x \Delta y \Delta z)(J_{\pmb{v}}^m)^TJ_{\pmb{r}}^m\pmb{x_r}, \label{eq:H12handle}\\
H_{21}\pmb{x_v} & \approx 2\beta(\Delta x \Delta y \Delta z)(J_{\pmb{r}}^m)^TJ_{\pmb{v}}^m\pmb{x_v}.
\end{align}
\end{subequations}
The motivation of using a function handle $H\pmb{x}$ instead of the matrix $H$ itself is to avoid numerical multiplication of two big matrices which could be time-consuming. With the function handle, for example, in \eqref{eq:H12handle} we can instead multiply a matrix by a vector twice to arrive at the desired computation.

As for the formulation of $J_{\pmb{v_j}}^k$ and $J_{\pmb{r_j}}^k$, by the structure of \eqref{eq:chain} we have that for $j=0,\cdots,m-1$ and $k>j$,
\begin{equation}\label{eq:Jv}
J_{\pmb{v_j}}^k = L^{-1}S(\pmb{v_{k-1}})R(\pmb{r_{k-1}})\cdots L^{-1}S(\pmb{v_{j+1}})R(\pmb{r_{j+1}})L^{-1}B(\pmb{\rho_j},\pmb{r_j})
\end{equation}
where
\begin{equation}
B(\pmb{\rho_j},\pmb{r_j}) = \frac{\partial}{\partial \pmb{v_j}}(S(\pmb{v_j})R(\pmb{r_j})\pmb{\rho_j})
\end{equation}
is dependent only on $\pmb{\rho_j}, \pmb{v_j}$ because $S(\pmb{v_j})$ linear to $\pmb{v_j}$, and 
\begin{equation}\label{eq:Jr}
\begin{aligned}[t]
J_{\pmb{r_j}}^k = &\Delta t\cdot L^{-1}S(\pmb{v_{k-1}})R(\pmb{r_{k-1}})\cdots L^{-1}S(\pmb{v_{j+1}})R(\pmb{r_{j+1}})L^{-1}S(\pmb{v_j})\text{diag}(\pmb{\rho_i}\odot\pmb{\chi_i}).
\end{aligned}
\end{equation}
With the explicit and recursive expressions in equations \eqref{eq:Jv} and \eqref{eq:Jr}, one can compute $J_{\pmb{v}}^m=[J_{\pmb{v_0}}^m|J_{\pmb{v_1}}^m|\cdots |J_{\pmb{v_{m-1}}}^m]$, $J_{\pmb{r}}^m=[J_{\pmb{r_0}}^m|J_{\pmb{r_1}}^m|\cdots |J_{\pmb{r_{m-1}}}^m]$ and their transpose multiplied with a vector in an iterative manner. $J_{\pmb{v}}^T$ and $J_{\pmb{r}}^T$ multiplied with a vector can also be computed recursively due to their lower-triangularity.

\subsection*{Algorithm}
With the analytic formulation of the gradient $g$ and the Hessian handle $H\pmb{x}$ given above, we can then utilize the Gauss-Newton method to find the optimal solution. See Algorithm~\ref{alg:algo} for the pseudo-code. 

If we have more than two successive given images, $\pmb{\rho_0^{img}},\pmb{\rho_1^{img}},\cdots,\pmb{\rho_{q-1}^{img}}$ which is inherent to DCE-MRI studies and corresponding given indicator functions $\pmb{\chi_0^{img}},\pmb{\chi_1^{img}},\cdots,\pmb{\chi_{q-2}^{img}}$ between adjacent images where $q>2$ and $q\in\mathbb{N}^+$, we can run the algorithm iteratively between each pair of adjacent images to derive prolonged velocity fields and relative sources. In other words, the process can be graphed as $\pmb{\rho_0^{img}} \xrightarrow[\text{loop 1}]{\text{urOMT}}\pmb{\rho_1^{img}} \xrightarrow[\text{loop 2}]{\text{urOMT}}\cdots\xrightarrow[\text{loop $q-2$}]{\text{urOMT}}\pmb{\rho_{q-2}^{img}} \xrightarrow[\text{loop $q-1$}]{\text{urOMT}}\pmb{\rho_{q-1}^{img}}$, where the urOMT algorithm is run for $q-1$ times and therefore $q-1$ successive outputs are returned. If one would like the prolonged velocity fields to be smoother in the temporal dimension, one can put the last interpolated image $\pmb{\rho_{m}}$ of the previous loop into the next loop as the initial image to avoid constantly introducing new data noise into the system\cite{chen2022visualizing1}.

\begin{algorithm}
\caption{Gauss-Newton Method}
\label{alg:algo}
\begin{algorithmic}[1]
\STATE{Load in $\pmb{\rho_0^{img}}, \pmb{\rho_1^{img}}$ and $\pmb{\chi}$, and other parameters;}
\STATE{$\pmb{x} = [\pmb{v}; \pmb{r}] =$ initial guess (all zeros);}
\FOR{$i = 1,2,\cdots, MaxIter$}
\STATE{Compute interpolations $\pmb{\rho}$ recursively from the discretized partial differential equation \eqref{eq:chain}: $\pmb{\rho}=$ SrcAdvDiff$(\pmb{\rho_0^{img}},\pmb{\chi},\pmb{v},\pmb{r})$;}
\STATE{Compute $S(\pmb{v_j})$, $R(\pmb{r_j})$ and $B(\pmb{\rho_j},\pmb{r_j})$ for $j=0,\cdots,m-1$;}
\STATE{Compute gradient $g$ and the Hessian function handle $H\pmb{x}$ according to equations \eqref{eq:gradients1}-\eqref{eq:hessians};}
\STATE{Solve linear system $H\pmb{x}=-g$ for $\pmb{x}$;}
\STATE{Do line search to find length $l$;}
\IF{line search fails}
    \RETURN $[\pmb{v};\pmb{r}]$;
\ENDIF
\STATE{Update $[\pmb{v};\pmb{r}]=[\pmb{v};\pmb{r}]+l\pmb{x}$;}
\ENDFOR
\RETURN $[\pmb{v};\pmb{r}]$;
\end{algorithmic}
\end{algorithm}

\section*{Results}
\label{sec:res}
In this section, we give some examples illustrating the application of the urOMT methodology. Obviously, the model admits both the Eulerian and Lagrangian perspectives for post-processing. The Eulerian formulation focuses on the fluid flows at fixed locations over time, while the Lagrangian formulation enables one to follow the trajectory of a given particle over time, and therefore to analyze the features along the given trajectory. 

Specifically, suppose we run the urOMT algorithm on given density images $\pmb{\rho_0^{img}}, \pmb{\rho_1^{img}},\cdots,\pmb{\rho_{q-1}^{img}}$ and given indicators $\pmb{\chi_0^{img}},\pmb{\chi_1^{img}},\cdots,\pmb{\chi_{q-2}^{img}}$ with discretization described before, the algorithm will run $q-1$ successive loops and return the solutions
\begin{equation}
\pmb{v_k^*} =
\begin{bmatrix}
\pmb{v_{k,0}^*}\\
\vdots\\
\pmb{v_{k,m-1}^*}
\end{bmatrix}
\text{and }
\pmb{r_k^*} =
\begin{bmatrix}
\pmb{r_{k,0}^*}\\
\vdots\\
\pmb{r_{k,m-1}^*}
\end{bmatrix}
\end{equation}
for $k=1,\cdots,q-1$ where the subscript $k$ stands for the $k$-th loop. By taking the $L^2$ norm of $\pmb{v_{k,0}^*},\cdots,\pmb{v_{k,m-1}^*}$, we derive the optimal speed
\begin{equation}
\pmb{s_k^*} =
\begin{bmatrix}
\pmb{s_{k,0}^*}\\
\vdots\\
\pmb{s_{k,m-1}^*}
\end{bmatrix}
\triangleq
\begin{bmatrix}
||\pmb{v_{k,0}^*}||\\
\vdots\\
||\pmb{v_{k,m-1}^*}||
\end{bmatrix}
.
\end{equation}
Therefore, $\pmb{r_{k,i}^*}$ and $\pmb{s_{k,i}^*}$ are called the \textit{Eulerian relative source maps} and \textit{Eulerian speed maps}, respectively, for $k=1,\cdots,q-1$ and $j=0,\cdots,m-1$, and they allow us to observe the fluid flows at fixed coordinates at different discrete time steps. For ease of visualization, we define the \textit{time-averaged Eulerian speed map} and \textit{time-averaged Eulerian relative source map} between $\pmb{\rho_{N_0}^{img}}$ and $\pmb{\rho_{N_1}^{img}}$ as
\begin{equation}
\frac{1}{m(N_1-N_0)}\sum_{k=N_0+1}^{N_1}\sum_{j=0}^{m-1} \pmb{s^*_{k,j}}
\text{\quad and \quad}
\frac{1}{m(N_1-N_0)}\sum_{k=N_0+1}^{N_1}\sum_{j=0}^{m-1} \pmb{r^*_{k,j}},
\end{equation}
respectively, where $0\leqslant N_0<N_1\leqslant q-1$ and $N_0,N_1\in\mathbb{N}^+$. 

In contrast, the post-processing framework developed and applied in the previous work\cite{koundal2020,chen2022,chen2022visualizing1,ozturk2023continuous} that follows Lagrangian coordinates present data as binary trajectories of the fluid flows, which we refer to as the \textit{pathlines} in our work. By connecting the starting and terminal points of the pathlines, we have the \textit{velocity flux vectors}, which are also called the \textit{displacement fields} in physics. These vectors may be used to visualize the direction and the distance travelled in a compact and interpretable manner. If we endow the pathlines with more information, for example, speed (the $L^2$ norm of the velocity field) and P\'{e}clet ($Pe$) number (i.e., the ratio of the rate of advection to diffusion), we can derive what we call \textit{speed-lines} and \textit{P\'{e}clet-lines}, respectively. More details of this Lagrangian method may be found in Koundal et al.\cite{koundal2020} and Chen et al.\cite{chen2022visualizing1}

Now that we have briefly described the two post-processing methods, the Eulerian and the Lagrangian perspectives, we next utilize two datasets to exhibit the results of the  urOMT analysis. The first dataset is a synthetic geometric dataset derived from Gaussian spheres as a simple demonstration of the urOMT method. The second is the DCE-MRI data from a rat brain where we elucidate the application of urOMT in \textit{in vivo} datasets for quantifying and visualizing the fluid flows. 

\subsection*{Tests on Gaussian Spheres}\label{sec:gauss}

\begin{figure}[htbp]
\begin{center}
\setlength{\tabcolsep}{1pt}
\includegraphics[width=0.95\textwidth]{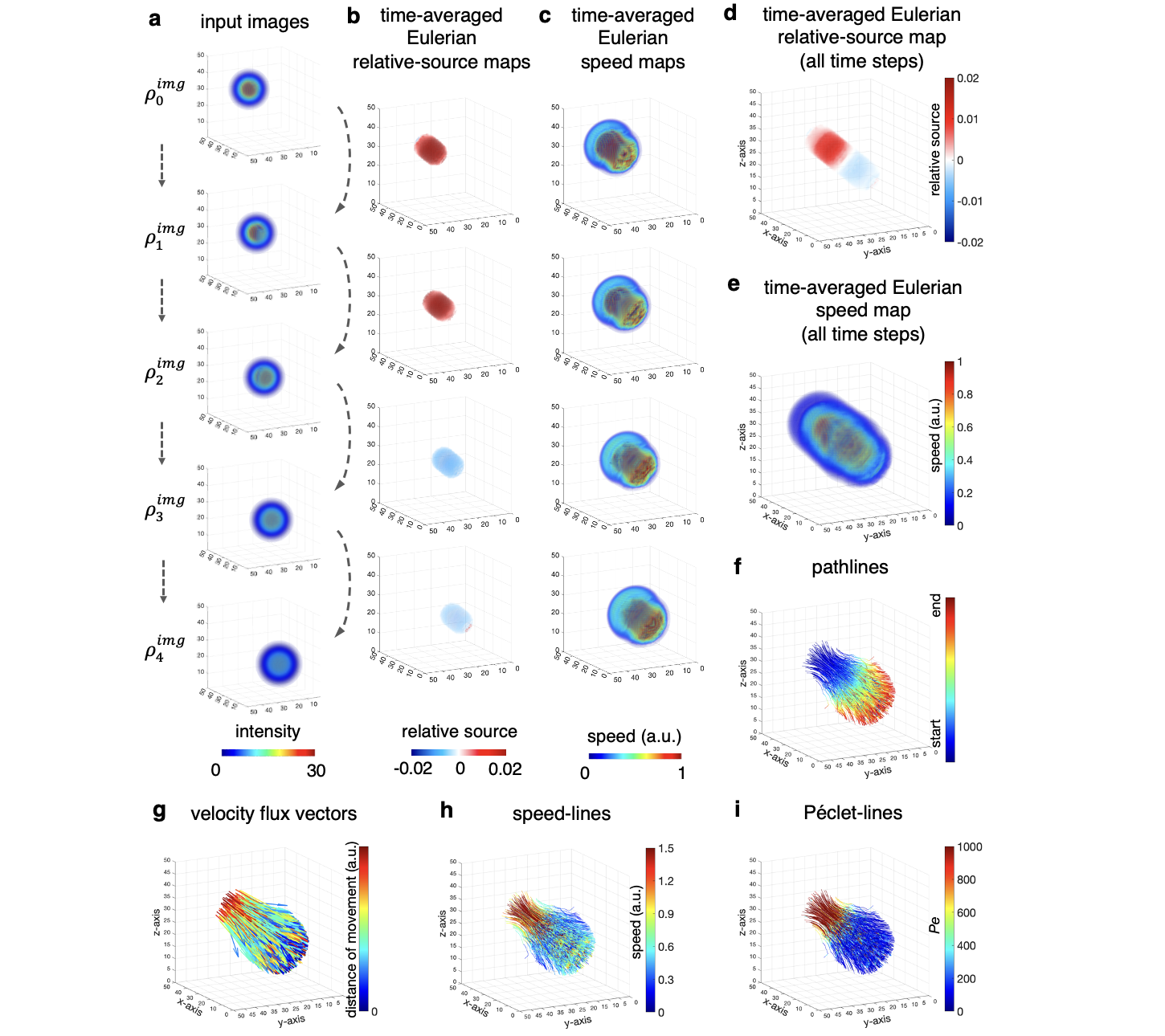}
\end{center}
\caption{The First Test on 3D Gaussian Spheres. \textbf{a}: Five successive images, shown in 3D rendering, were created from Gaussian spheres as inputs into the urOMT algorithm. In addition to advection (from top-left to bottom-right) and diffusion included in the transport process, mass was gained from $\rho_0^{img}$ to $\rho_2^{img}$ and was lost from $\rho_2^{img}$ to $\rho_4^{img}$ in the center region. \textbf{b-e}: Under Eulerian coordinates, the time-averaged speed maps and relative source maps visualized in 3D indicate the speed and mass gain/loss distribution in the domain, respectively. \textbf{b} and \textbf{c} are derived between every pair of input images; \textbf{d} and \textbf{e} are derived between $\rho_0^{img}$ and $\rho_4^{img}$. \textbf{f-i}: Under Lagrangian coordinates, the binary trajectories of the transport are recorded by pathlines color-coded with start and end points. Connecting the start and end points of pathlines, we derive the velocity flux vectors illustrating the direction and distance of the overall movement. By endowing the pathlines with speed and P\'{e}clet ($Pe$) number, we derive the speed-lines and P\'{e}clet-lines, respectively.}
\label{fig:result_gauss}
\end{figure}

We first start with creating five successive 3D images from Gaussian spheres denoted as $\rho_0^{img},\rho_1^{img},\cdots,\rho_4^{img}$ in order to test the urOMT algorithm. These spheres were created to first gain mass and later lose mass in the center region of the spheres over time, and in addition the spheres are spatially transported forward and are also exhibiting active diffusion over time (Figure~\ref{fig:result_gauss}a). Specifically, $\rho_i$ was created from a 3D Gaussian function

\begin{equation}
   G_i(x,y,z) = \frac{100}{\sqrt{2\pi}}\exp{\big(-\frac{(x-0.8i)^2+(y-0.8i)^2+(z-0.8i)^2}{2}\big)},
\end{equation}
then advection is naturally included in the whole process to reflect the forward translation of the center of the Gaussian spheres from top-left to bottom-right (Figure~\ref{fig:result_gauss}a). We define the center region of these Gaussian spheres as the region within a radius of 1.5 of the center $[0.8i,0.8i,0.8i]$ which are denoted as $\chi_i(x,y,z)$ for $i=0,\cdots,4$. Then we apply $\chi_i$ to $G_i$ to imitate mass gain and loss in the center region which are given by $(1+a_i\chi_i(x,y,z))G_i(x,y,z)$ where $[a_0,a_1,a_2,a_3,a_4] = [0,0.1,0.2,0.1,0]$. We then discretize these functions by taking a uniform spatial length and fitting them into the same numerical grid of size $50\times50\times50$. Diffusion was further added to $\rho_i$ by applying a MATLAB inbuilt 3D Gaussian filter \textit{imgaussfilt3} to $\rho_i$ with standard deviation = $(i+1)\sqrt{0.2}$ for $i = 1,\cdots,4$. 

One can imagine the five successive images as visualizing a ``wormhole'', which is moving forward (advection) with mass diffusing into its surrounding area. At the same time, at the center region of the wormhole mass is gained ($\rho_0^{img}$ to $\rho_2^{img}$) or lost ($\rho_2^{img}$ to $\rho_4^{img}$) from or to another interconnected space.

In our first test on the Gaussian sphere data, we fed all five images, $\rho_0^{img},\rho_1^{img},\cdots,\rho_4^{img}$, into the urOMT algorithm. From $\rho_i^{img}$ to $\rho_{i+1}^{img}$, we utilize the indicators as the linear translation of $\chi_{i}$ to $\chi_{i+1}$ (the region within a radius of 1.5 centered at $[0.8(i+\frac{j}{m}),0.8(i+\frac{j}{m}),0.8(i+\frac{j}{m})]$ for $j=0,\dots,m, i=0,\dots,3$) to only allow mass gain and loss to occur in the center regions. The parameters used in this experiment are listed in Table~\ref{tab:param}. Computations were run with MATLAB 2018b on the departmental High Performance Computing cluster at Memorial Sloan Kettering Cancer Center with Red Hat Enterprise Linux 7.5 operating system using 3 CPUs and 128GB of memory which took about 3 hours and 20 minutes.

From the Eulerian perspective, we show the time-averaged Eulerian speed maps and relative source maps between every pair of input images (Figure \ref{fig:result_gauss}b,c) and between $\rho_0^{img}$ and $\rho_4^{img}$ (Figure \ref{fig:result_gauss}d,e). The algorithm recognized the core regions of the spheres as having higher speed. Moreover, it successfully captured the mass gain and loss patterns during the entire process with the red color indicating the initial mass gain from $\rho_0^{img}$ to $\rho_2^{img}$ and the blue color indicating the later mass loss from $\rho_2^{img}$ to $\rho_4^{img}$. The Eulerian relative source map is restricted in the center region of the spheres because of the indicator functions. With Lagrangian post-processing, we derived the pathlines, velocity flux vectors, speed-lines and P\'{e}clet-lines under Lagrangian coordinates (Figure~\ref{fig:result_gauss}f-i). The pathlines exhibit what trajectories of particles would look like over time if they were placed at given initial points in the system at $t=0$. The funnel-shape of the pathlines is a reflection of the accumulated effect of diffusion which gradually disperses the mass. The velocity flux vectors are also provided to show the direction and distance of the whole transport process. The speed-lines, i.e., pathlines endowed with speed, indicate that higher speed occurred mainly in the core of the spheres which is in agreement with the Eulerian speed map. The P\'{e}clet-lines show that initially the transport was dominated by advection and later by diffusion.

\begin{figure}[htbp]
\begin{center}
\setlength{\tabcolsep}{1pt}
\includegraphics[width=0.95\textwidth]{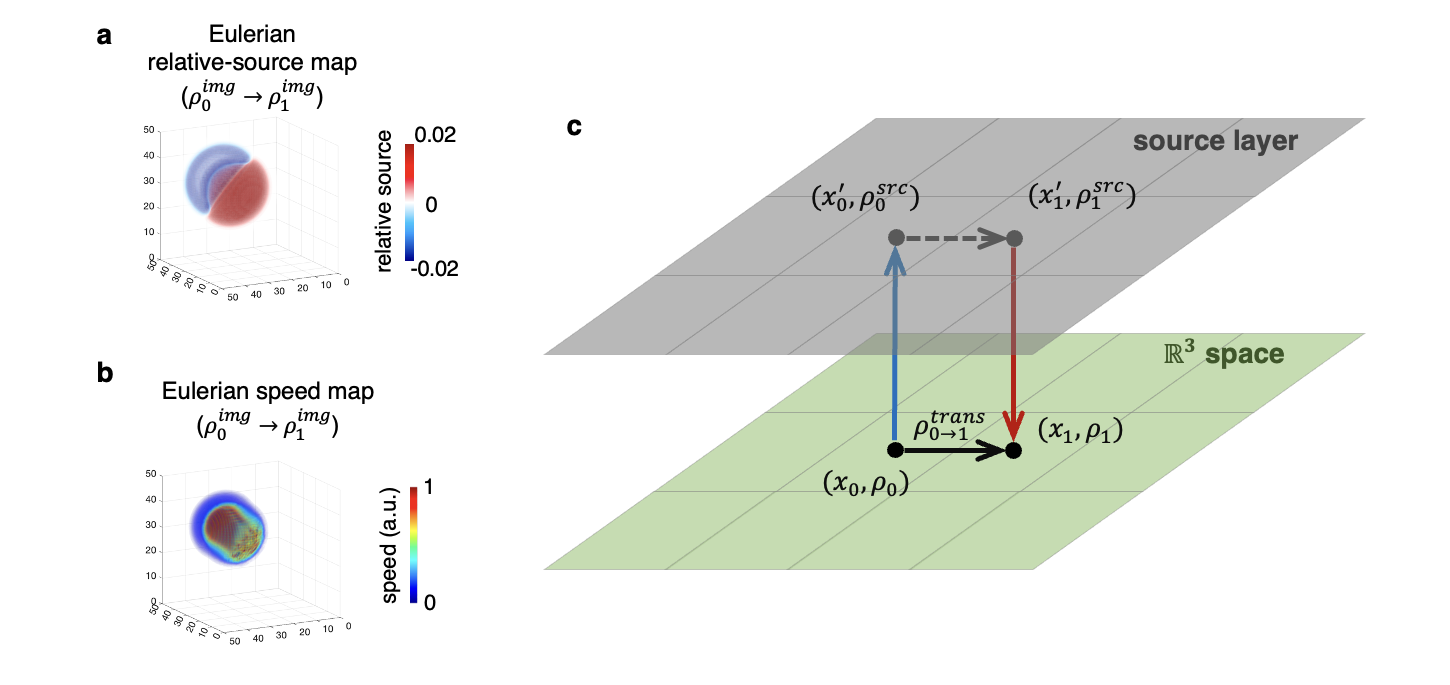}
\end{center}
\caption{The Second Test on 3D Gaussian Spheres. \textbf{a-b}: Given input images $\rho_0^{img}$ and $\rho_1^{img}$ visualized in Figure~\ref{fig:result_gauss}a, Eulerian outputs, shown in 3D rendering, were returned from the urOMT analysis. \textbf{c}: The transport in the system can be separated into two channels where the advection and diffusion take place in the $\mathbb{R}^3$ space and the relative source pushes or draws mass between the two channels.}
\label{fig:result_gauss2}
\end{figure}

To further understand how indicators affect the results and how mass is transferred in the urOMT system, we performed a second test with only two input images $\rho_0^{img},\rho_1^{img}$ and indicators as a constant 1 to allow mass gain/loss everywhere. The parameters used in this test are listed in Table~\ref{tab:param}. The computational runtime was about 30 minutes with the same machine and configuration as the first test. Eulerian relative source map (Figure~\ref{fig:result_gauss2}a) and Eulerian speed map (Figure~\ref{fig:result_gauss2}b) were returned from the urOMT algorithm. We observe that the top-left of the Eulerian relative source map was negative and colored in blue where the input sphere leaves and the bottom-right was positive colored in red where the sphere arrives. In other words, without any constraint on the relative source term, the Eulerian relative source map can provide global information on the transport status in terms of mass loss and arrival. 

In this case, mass is freely transferred via two ``channels'', the real $\mathbb{R}^3$ space and the imaginary source layer (Figure~\ref{fig:result_gauss2}c). Advection and diffusion occur only in the $\mathbb{R}^3$ space. The source layer can generate infinite amount of mass and push it to the corresponding location in the $\mathbb{R}^3$ space; it can also draw a certain amount of the mass from the $\mathbb{R}^3$ space to the corresponding location in the source layer. The activities in both channels happen simultaneously. Suppose in $\mathbb{R}^3$, one would like to move mass $\rho_0>0$ at location $x_0$ to mass $\rho_1>0$ at location $x_1$, i.e., 
$\begin{bmatrix}
(x_0,\rho_0)\\ (x_1,0)
\end{bmatrix}
\Rightarrow
\begin{bmatrix}
(x_0,0) \\ (x_1,\rho_1)
\end{bmatrix}$
.
We denote $x_0^\prime$ and $x_1^\prime$ as the corresponding locations of $x_0$ and $x_1$ in the source payer, respectively; the amount of mass transported in $\mathbb{R}^3$ from $x_0$ to $x_1$ as $\rho_{0\rightarrow1}^{trans}>0$; the amount of mass drawn from $x_0$ to $x_0^\prime$ as $\rho_0^{src}>0$; and the amount of mass pushed from $x_1^\prime$ to $x_1$ as $\rho_1^{src}>0$. Therefore, the following equations hold
\begin{equation}
\rho_0 = \rho_{0\rightarrow1}^{trans} + \rho_0^{src}, \quad
\rho_1 = \rho_{0\rightarrow1}^{trans} + \rho_1^{src}.
\end{equation}
The weighting parameter $\alpha$ in the urOMT formulation therefore balances the split of $\rho_0$ into $\rho_{0\rightarrow1}^{trans}$ and $\rho_0^{src}$ whose effect will be demonstrated in the next section. 

In this test, to move $\rho_0^{img}$ to $\rho_1^{img}$, part of the mass was transported in $\mathbb{R}^3$ and part in the source layer, and these two phenomena facilitate each other. We know \textit{a priori} from the creation of the data that the sphere should be transported forward and mass gain only occurs in the center region. However, allowing a global relative source, the urOMT algorithm finds an easier way to transform $\rho_0^{img}$ into $\rho_1^{img}$ by making use of the source layer to transfer mass.

Most importantly, removing the indicator constraint the relative source globally compensates the transport in $\mathbb{R}^3$ and is indicative of the leaving and arrival of mass. This can be very useful in the context of applications to some real dataset. As a matter of fact, \textit{a priori} indicator for the source layer is sometimes unavailable due to the complexity of the system.

\begin{table}[htbp]
\begin{center}
    \begin{tabular}{ ||m{1.8cm}|m{4.7cm}|m{2.8cm}|m{2.8cm}|m{2.8cm}|  }
    \cline{1-5}
    \textbf{Parameter} & \textbf{Definition} & \textbf{Value for Gaussian Sphere Data Test 1} & \textbf{Value for Gaussian Sphere Data Test 2} & \textbf{Value for Brain Data}      \\ \cline{1-5}
    \hline\hline
    $n_1$ & grid size in $x$ axis & \multicolumn{2}{c|}{50} & 56 \\ \cline{1-5}
    $n_2$ & grid size in $y$ axis & \multicolumn{2}{c|}{50} & 106 \\ \cline{1-5}
    $n_3$ & grid size in $z$ axis & \multicolumn{2}{c|}{50} & 51 \\ \cline{1-5}
    $q$ & number of input images & 5 & 2 & 15 \\ \cline{1-5}
    $m$ & number of time intervals between two input images& \multicolumn{3}{c|}{10} \\ \cline{1-5}
    $\Delta t$ & temporal spacing  & \multicolumn{3}{c|}{0.4} \\ \cline{1-5}
    $\Delta x$ & $x$-axis spacing & \multicolumn{3}{c|}{1}\\ \cline{1-5}
    $\Delta y$ & $y$-axis spacing & \multicolumn{3}{c|}{1}\\ \cline{1-5}
    $\Delta z$ & $z$-axis spacing & \multicolumn{3}{c|}{1}\\ \cline{1-5}
    $\sigma$ & diffusion coefficient & \multicolumn{3}{c|}{0.002} \\ \cline{1-5}
    $\alpha$ &  weighting parameter for the source term & 9000 & \multicolumn{2}{c|}{10000} \\\cline{1-5}
    $\beta$ &  weighting parameter for the fitting term & \multicolumn{2}{c|}{5000} & 50 \\\cline{1-5}
    $\chi$ &  indicator function of the relative source $r$ & 1's in the center regions, otherwise 0's & \multicolumn{2}{c|}{all 1's} \\\cline{1-5}
    \end{tabular}
\caption{Parameters used in the urOMT algorithm.}
\label{tab:param}
\end{center}
\end{table}

\subsection*{Application to Rat Brain MRI}\label{sec:rat}
A very useful application of the urOMT model is to quantify the transport properties of brain fluids with DCE-MRI data. It still remains debated in the scientific community how fluid and solutes are transported in brain parenchyma and how the "dirty" fluid is drained out of the brain to maintain homeostasis\cite{BOHR2022104987,zhao2022physiology}. The relative source in the urOMT model may be helpful in revealing fluid and solute clearance patterns in the brain.

In this experiment, our urOMT method was applied to 3D DCE-MRI dataset derived from a 3-month-old healthy rat brain. The tracer, gadoteric acid, was injected into the cerebrospinal fluid (CSF) of the rat after the rat was anesthetized. The DCE-MRI data series of the rat brain was collected every 5 minutes and lasted for 140 minutes, ending up with 29 images in total. The MRI signal images were then processed to derive the \%-signal change from the baseline to approximate the concentration images of tracers. More information of the DCE-MRI data may be found in Chen et al.\cite{chen2022} 

\begin{figure*}[htbp]
\begin{center}
\includegraphics[width=0.95\textwidth]{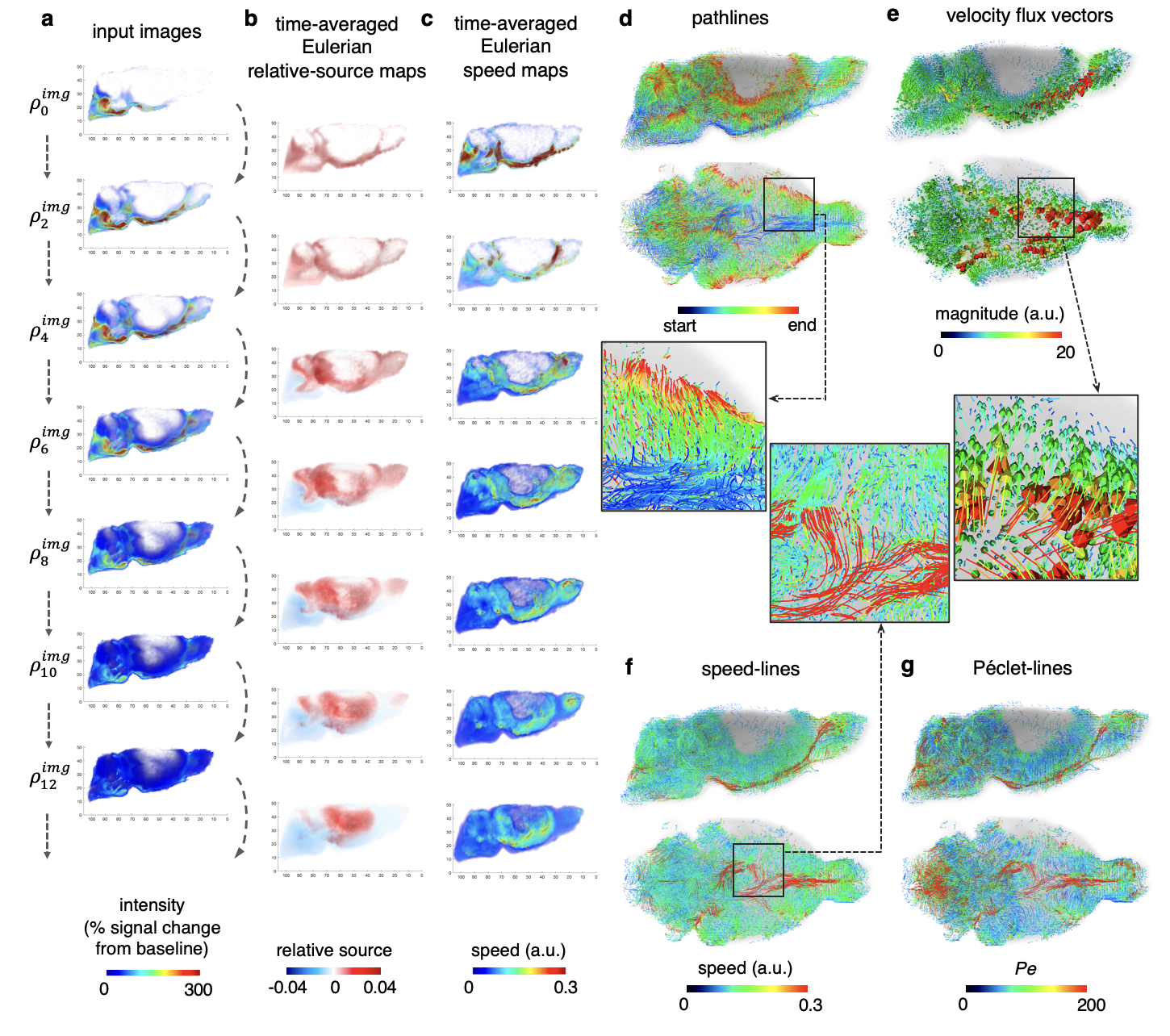}
\end{center}
\caption{Application to 3D Rat Brain MRI. \textbf{a}: Rat brain MRIs, shown in 3D rendering from the right-lateral view plane, were fed successively into the urOMT algorithm. \textbf{b-c}: As returned outputs, the Eulerian time-averaged relative source maps and speed maps between every other image were plotted, indicating the rate of mass gain/loss and speed distribution over time, respectively. \textbf{d-g}: Under Lagrangian coordinates, pathlines, color-coded with start and end points, show the trajectories of the tracers in brain. The speed-lines show the speed values along pathlines, and similarly P\'{e}clet-lines indicate whether the transport is advection or diffusion-dominated along pathlines. The velocity flux vectors show the direction and distance of the transport. All lines and vectors are shown from both the right-lateral and bottom views, and are overlaid on the anatomical data in gray.}
\label{fig:result_rat}
\end{figure*}

\begin{figure*}[htbp]
\begin{center}
\includegraphics[width=0.95\textwidth]{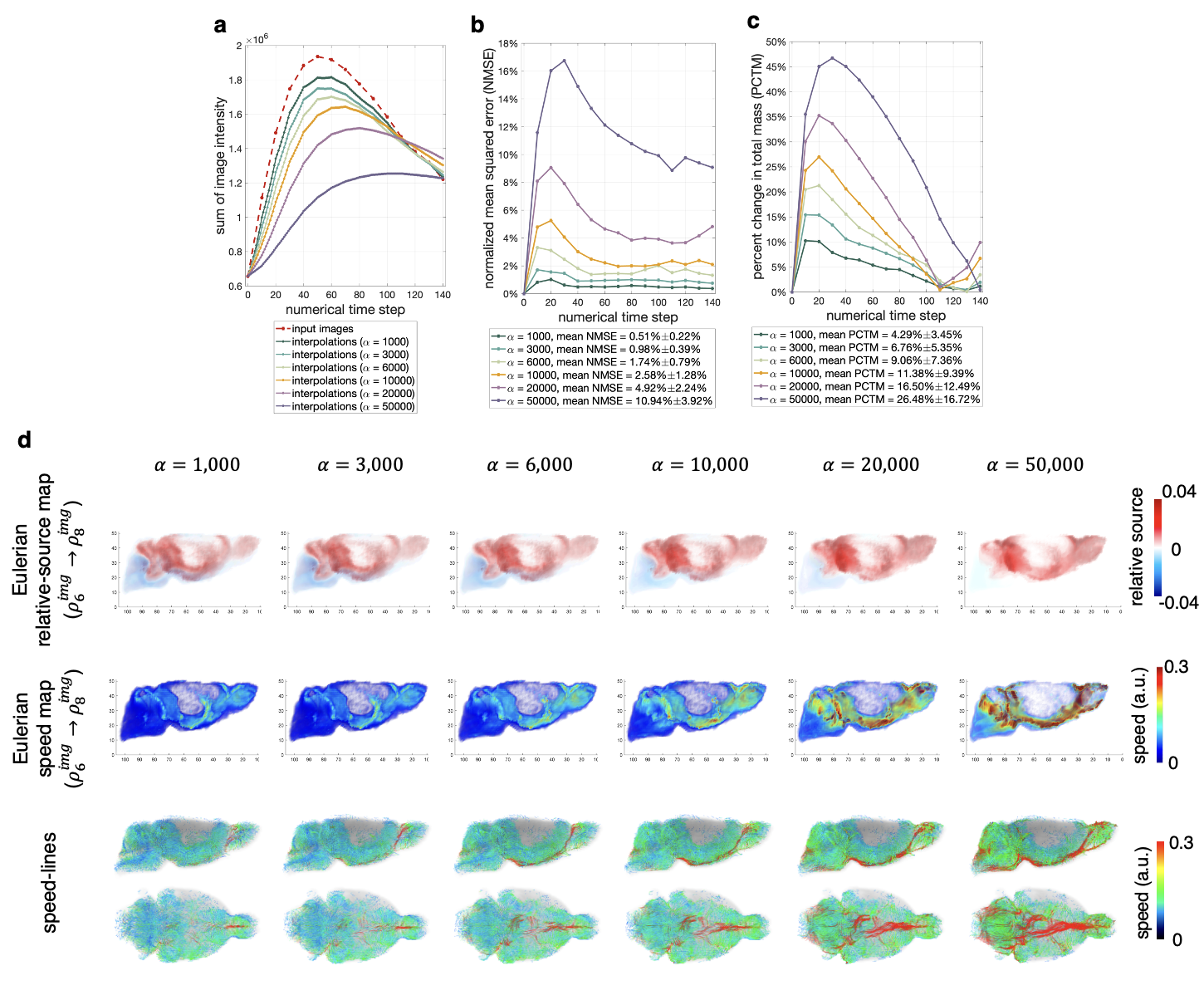}
\end{center}
\caption{Examination of the Effect of the Parameter $\alpha$. \textbf{a}: The comparison of the total image intensity curve of the input images and interpolations with different $\alpha$'s. \textbf{b}: The normalized mean squared error curve over the numerical steps with different $\alpha$'s. \textbf{c}: The percent change in total mass curve over the numerical steps with different $\alpha$'s. \textbf{d}: The comparison of the Eulerian relative source maps (first row) and Eulerian speed maps (second row) from $\rho_6^{img}$ to $\rho_8^{img}$ (shown in Figure~\ref{fig:result_rat}a) and the speed-lines (the third row) returned from the urOMT and its post-processing algorithm with different $\alpha$'s. The first two rows are from the right-lateral view and the third row is from both the right-lateral and bottom views.}
\label{fig:performance}
\end{figure*}

In this numerical experiment, we assume that the intensity of the DCE-MRI images is proportional to the mass in the urOMT model. We used every other image within the brain region as input images in order to save running time, which resulted in a total of 15 images: $\rho_0^{img},\rho_1^{img},\cdots,\rho_{14}^{img}$ for the urOMT algorithm (Figure~\ref{fig:result_rat}a). In the previous work\cite{chen2022}, the rOMT model was applied right before the peak of the total signal intensity of the input images, i.e., $\rho_3^{img},\cdots,\rho_{14}^{img}$ in current notation. In the present experiment, due to the introduction of the relative source term in the model we were able to include earlier frames when the total intensity was still rapidly increasing (Figure \ref{fig:performance}a, the red dashed curve). Since the mechanism of the fluid transport in brains is complex and still under intense investigation, there is no information provided \textit{a priori} for the relative source term. So we set its indicator function $\chi$ to be equal to 1 everywhere with the assumption that the entire brain system is "leaky" to allow the tracer to enter and exit the system through unknown ways. In order to derive smooth prolonged dynamics, we used the last interpolated image from the previous numerical loop as the initial image in the next loop. The parameters used in this experiment are listed in Table \ref{tab:param}. The computation was also run with MATLAB 2018b on our departmental High Performance Computing cluster at Memorial Sloan Kettering Cancer Center with Red Hat Enterprise Linux 7.5 operating system using 3 CPUs and 128GB of memory which took about 8 hours and 50 minutes.

We show the Eulerian relative source maps and speed maps between every other input frame (Figure~\ref{fig:result_rat}b-c). The Eulerian relative source maps indicate that the tracer first entered the CSF surrounding the brain causing the intensive mass gain noted in the first frames, and then subsequently moved into deeper brain tissue regions resulting in mass loss in CSF and mass gain in the brain tissue. The Eulerian speed maps indicate that the initial speed of tracer when entering the CSF was very high, and slowed when the tracer penetrated deeper, probably due to diffusion dominating the transport in the tissue. From the Lagrangian perspective, we show the pathlines of tracers starting at $t=0$ as well as those lines endowed with speed and P\'{e}clet number (Figure \ref{fig:result_rat}d,f-g). The pathlines signifies the pathways/trajectories of the tracer entering CSF and brain tissue. The speed-lines show that higher speed was mainly in the CSF and along the perivascular space of the large vessels, and further that the transport was identified as advection-dominated according to the high values in the P\'{e}clet-lines in those regions. The velocity flux vectors were also derived and demonstrated that the tracer entered the brain tissue in a symmetrical pattern about the the midline of the skull base  (Figure \ref{fig:result_rat}e).

Recall that in the urOMT model \eqref{eq:uromt_dp1}-\eqref{eq:uromt_dp3}, the weighting parameter $\alpha>0$ penalizes the source term in the cost function. In theory, as $\alpha \rightarrow +\infty$, $r$ gets suppressed and this model approximates the rOMT model where unbalanced mass gain and loss are not allowed. In the test, we used above $\alpha=10000$. We ran further tests with different values including $\alpha=1000, 3000, 6000, 20000$ and $50000$ to demonstrate the effect of $\alpha$ on the results.

From the urOMT algorithm, we can compute the final interpolation $\rho_{i,m}$ for $i=1,\cdots,q-1$.
To test the numercial performance, we define the \textit{normalized mean squared error} (NMSE) between each pair of final interpolation $\rho_{i,m}$ and the ground truth $\rho_{i}^{img}$ as
\begin{equation}
\frac{ \| \rho_{i,m}-\rho_{i}^{img} \|^2}{\|\rho_{i}^{img}\|^2}\times100\%.
\end{equation}
We also define the \textit{percent change in total mass} (PCTM) between $\rho_{i,m}$ and $\rho_{i}^{img}$ as
\begin{equation}
\frac{ \lvert \text{sum}(\rho_{i,m})-\text{sum}(\rho_{i}^{img}) \rvert}{\text{sum}(\rho_{i}^{img})}\times100\%
\end{equation}
where $\text{sum}(\cdot)$ denotes taking the sum of all entries in a vector. Both NMSE and PCTM measure the accuracy of the urOMT algorithm and the lower they are, the more accurate the urOMT results are. However, NMSE measures the closeness of the ground truth and the interpolations in a local manner, while PCTM in a global manner. We emphasize that urOMT is a data-driven method; in other words, the transport between two images is a black box and urOMT is simply giving the most likely evolving path under a pre-defined least cost assumption. Therefore, it is difficult to define the ground truth between two input images for us to compare with.

By fixing the rest of the parameters but using various $\alpha$ values (1000, 3000, 6000, 10000, 20000 and 50000), we ran our novel urOMT algorithm and the post-processing procedure on the same rat brain data. In Figure~\ref{fig:performance}a, we plot the curves of the total image intensity of both input images and interpolations with different $\alpha$'s. From the curves, the smaller the $\alpha$ is, the closer the curve of interpolations is to the curve of input images. This is in agreement with the theoretical observation that the smaller the $\alpha$ is, the less the source term is penalized, which means more mass gain/loss is allowed in the system in order to adjust the total mass. From Figure~\ref{fig:performance}b-c, we found that a smaller $\alpha$ may help derive more accurate numerical simulations, because for both NMSE and PCTM the trend is that the larger the $\alpha$ is, the higher the NMSE and PCTM values are. This makes sense because when the system is highly unbalanced, and forcefully using a high $\alpha$ produces a nearly balanced environment which contradicts with the data setting. To further examine the effect of $\alpha$ on the quantitative and visual results, we plot the Eulerian relative source maps and Eulerian speed maps from $\rho_6^{img}$ to $\rho_8^{img}$ and speed-lines with different $\alpha$'s in Figure \ref{fig:performance}d. They show that when the source term is greatly penalized with a high $\alpha$, the speed is consequently elevated. For example in Figure~\ref{fig:performance}d with $\alpha = 50000$ where the model most approximates rOMT, there is very high speed transport color coded in red at the base of the brain, which could be artificial and potentially over-estimated because in such a system mass has to move more quickly in order to match the final input image. In contrast, in a system with low $\alpha$ where instantaneous mass gain/loss is promoted, mass does not need to transport in the same degree to match the final input image, since it can instead pull in (or push out) mass from (or into) the ``invisible sink" (the source layer in Figure~\ref{fig:result_gauss2}c) via the relative source $r$. For example in Figure~\ref{fig:performance}d with $\alpha = 1000$, the rapid movement of tracer at the base of CSF almost disappeared because the change of the system is mostly accounted by the relative source. 
In general, the relative source $r$ and the velocity field $v$ compensate each other and their effects in the system are actively controlled by the weighting parameter $\alpha$.


\section*{Discussion}\label{sec:discussion}
In this work, we introduced an unbalanced version of the rOMT model for studying brain fluid dynamics using DCE-MRI images, which we referred to as urOMT. This method was utilized to make rOMT \cite{koundal2020,chen2022,chen2022visualizing1} more physically and biologically relevant by removing the mass conservation constraint. Specifically, the urOMT model accounts for the change of the total mass in the system by adding an independent relative source term into the formulation, while the rOMT model requires the total mass to be conserved. As discussed before, both theoretically and numerically the rOMT model can be approximated by the urOMT model when the parameter $\alpha$ goes to infinity. In other words, urOMT ``incorporates'' rOMT, and one can use urOMT with a large enough $\alpha$ (thereby assuring that the mass conservation condition of input images to be met) in place of the previous version of rOMT. As such the new urOMT introduced in the present work is a more powerful and flexible model for analysis of fluid transport.

The urOMT method may be particularly useful for studying the cross-talk between the glymphatic system and meningeal lymphatics. See relevant work\cite{hablitz2021glymphatic,BOHR2022104987,benveniste2019glymphatic} for more details about how the two systems interplay. With the additional information of the relative source in our urOMT model which reveals mass gain and loss, we are now able to observe and quantify the solute and fluid entering and exiting the two systems. For example, in Figure~\ref{fig:result_rat}b, at first we see mainly mass gain (colored in red) indicating that the tracer is flowing into the CSF, and then slowly we see blue color along the skull base, indicating that the tracer is either redistributing into the tissue bed or exiting via the draining lymphatic vessels. Therefore, our urOMT method has the potential to probe the clearance pattern in more depth at the level of the CSF and tissue compartments.

In the test on the rat brain DCE-MRI, we posed no spatial constraint on the relative source $r$ and used an indicator $\chi$ all 1's given that so far there is no agreed upon answer yet on the underlying transport mechanisms and the exact drainage pathways from the brain \cite{BOHR2022104987,zhao2022physiology}. Indeed from the DCE-MRI data (Figure~\ref{fig:result_rat}a), we did not observe a specific anatomical efflux route of the tracer out of the brain but the total intensity curve did decline in later frames (Figure~\ref{fig:performance}a, the red dashed curve). In this case, we assume that the system is leaky and the tracer can be transferred via unknown tunnels in brain (a correspondence of the source layer in Figure~\ref{fig:result_gauss2}c) to form the given images. Similarly in the widely used Toft's model\cite{tofts1997modeling,tofts1999estimating} for pharmacokinetic analysis of DCE-MRI studies of tumors where the signal curves are also unbalanced, a leakage between the blood vessels and the tissue is assumed to occur everywhere in the images. The popular parameter $K^{trans}$ returned from the Toft's model is used to quantify the local leakage of gadolinium-based tracers from the blood to the tissue\cite{tofts1997modeling,tofts1999estimating}.

In the urOMT algorithm, the difference between two input images is mainly captured by either the mass gain/loss rate $r$ in the source layer or the velocity field $v$ in $\mathbb{R}^3$, and $\alpha$ is the parameter balancing the two. Indeed as demonstrated by other work\cite{chizat2018scaling}, with a decreasing $\alpha$ parameter the transport (characterized by the velocity field $v$) is being compensated by mass gain and loss (characterized by the relative source $r$). With the rat brain MRI dataset, we demonstrated that a low $\alpha$ value gives higher numerical accuracy, but produces decreased speed by allowing the relative source term to play a greater role in the dynamics. Thus, one needs to be aware of the trade-off between the accuracy and the strength of fluid flows by choosing an appropriate parameter $\alpha$ when applying urOMT. One approach would be to make use of the indicator function $\chi$ to restrict the behaviors of the relative source $r$ within a certain region if the entering and exiting information of the fluid is known beforehand. Another approach is to use a time-varying $\alpha$ given that the total intensity curve of the input image is usually not linearly increasing or decreasing. Indeed, we plan to explore this possibility in some future work.

Other than the parameter $\alpha$, there are also additional parameters worthy of examining, such as the weighting parameter $\beta$ for the fitting term in the cost function \eqref{eq:uromt2_func} and the constant diffusion coefficient $\sigma$. Some preliminary efforts have been made in using a non-linear diffusion term in the partial differential equation \eqref{eq:Diff}, given the non-constant diffusion phenomena in brain fluid flows \cite{xu2023}. Future work includes adding a non-linear and spatially dependent diffusion term in the current urOMT formulation.

The urOMT model was largely motivated by the changing pattern of the total signal curve from DCE-MRI experiments. Specifically in rat brains, the temporal signal intensity curve typically peaks at approximately 1 hour after tracer injection, and later on, either keeps decreasing (Figure \ref{fig:performance}a, red dashed curve) or reaches a plateau (this depends on the total amount of tracer administered), signifying that the total mass may be highly unbalanced over time in the system. Given that the DCE-MRI protocol has been widely used in cancer imaging in clinics \cite{yankeelov2007dynamic,turkbey2010role,hylton2006dynamic}, the urOMT method also has the potential to be applied to tumor DCE-MRI  data in human experiments to investigate the tumor vasculature, and to help pave the way for new medical treatments.


\section*{Conclusions}\label{sec:conclusions}
The urOMT methodology incorporates both advection and diffusion motions into the transport process, as well as allowing for mass gain and loss in the dynamic images by introducing a relative source variable. For special cases, it may also constrain the relative source to a given region or time interval. As an extension of the rOMT model \cite{koundal2020,chen2022,chen2022visualizing1}, the urOMT model removes the total mass conservation constraint, while keeping the attractive advection-diffusion framework, making it applicable to modeling the fluid flows in the brain under DCE-MRI protocol and many other real-world modeling problems in computational fluid dynamics.







 





\section*{Data Availability}
The code for the urOMT algorithm and its post-processings is available at \href{https://github.com/xinan-nancy-chen/urOMT}{https://github.com/xinan-nancy-chen/urOMT}. The synthetic data and the DCE-MRI data used in this work is also provided therein.

\bibliography{references}



\section*{Acknowledgements}
This research was supported by AFOSR grants (FA9550-20-1-0029,FA9550-23-1-0096), Army Research Office grant (W911NF2210292), NIH grant (R01AT011419), and a grant from the Cure Alzheimer's Foundation.

\section*{Author Contributions Statement}
The numerical method and code development were done by X. C. The theory of OMT and other material preparations were contributed by A. R. T. The DCE-MRI data collection and processing were performed by H. B. The first draft of the manuscript was written by X. C. and all authors commented on previous versions of the manuscript. All authors read and approved the final manuscript.


\section*{Additional Information}

\textbf{Competing Interests Statement}:
The authors declare that they have no competing interests.






\end{document}